\documentclass[11pt]{article}

\usepackage{amssymb,amsmath,color}
\usepackage{graphicx,pstricks}

\newtheorem{theo}{\bfseries \hs Theorem}[section]

\newtheorem{prop}[theo]{\bfseries \hs Proposition}

\newtheorem{lemma}[theo]{\bfseries \hs Lemma}
\newtheorem{corol}[theo]{\bfseries \hs Corollary}
\newtheorem{con}[theo]{\bfseries \hs Conjecture}

\numberwithin{equation}{section} % Automatically number equations within sections

\def\C{\mathord{\mathbb C}}
\def\R{\mathord{\mathbb R}}

\def\Z{\mathord{\mathbb Z}}
\def\N{\mathord{\mathbb N}}

  \newcommand{\e}{\mathbf{e}}

  \newcommand{\w}{\mathbf{w}}
  
  \newcommand{\x}{\mathbf{x}}
  
  \newcommand{\y}{\mathbf{y}}
  
  \newcommand{\z}{\mathbf{z}}
  \newcommand{\0}{\mathbf{0}}
  
  \newcommand{\trans}{^\top}
  \newcommand{\rank}{\mathrm{rank\;}}
  \newcommand\supp{\mathop\mathrm{supp}\nolimits}
  \newcommand\eq{\mathrel\Leftrightarrow}

\def\hs{\hspace*{\parindent}}
\def\proof{{\hs \bf Proof.\ }}
\def\qed{\hfill$\Box$}

 \title{Upper bounds on the magnitude of solutions of certain linear systems with integer coefficients}
 \author{Pedro J.\ Freitas\\
 Centro de Estruturas Lineares e Combinat\'oria\\
 Departamento de Matem\'atica -- FCUL,
 Universidade de Lisboa\\
 \texttt{pedro@ptmat.fc.ul.pt}\medskip\\
 Shmuel Friedland\\
 Department of Mathematics, Statistics, and Computer Science\\
 University of Illinois at Chicago\\
 Chicago, Illinois 60607-7045, USA\\
 \texttt{friedlan@uic.edu}\medskip \\
 Gaspar Porta\\
 Washburn University\\
 Topeka, Kansas\\
 \texttt{gaspar.porta@washburn.edu}
 }
 \date{May 4, 2012}
 
\begin{document}

 \maketitle
 \begin{abstract}
 In this paper we consider a linear homogeneous system of $m$ equations in $n$ unknowns with integer coefficients over the reals.
 Assume that the sum of the absolute values of the coefficients of each equation does not exceed $k+1$ for some positive integer
 $k$.  We show that if the system has a nontrivial solution then there exists a nontrivial solution $\x=(x_1,\ldots,x_n)\trans$
 such that $\frac{|x_j|}{|x_i|}\le k^{n-1}$ for each $i,j$ satisfying $x_ix_j\ne 0$.  This inequality is sharp.

 We also prove a conjecture of A.\ Tyszka related to our results.
 \end{abstract}


 \section{Introduction}

% Diophantine problems, as Archimedes' cattle problem \cite{Len02}, are of great mathematical interest for thousands of years.
% Usually, the classical problems are stated with a prescribed number of variables and equations, which are usually quite small.

 In this paper we consider $m$ homogeneous linear equations, with integer coefficients, in $n$ variables.  I.e., our system is $A\x=\0$, where $A$ is an $m \times n$ matrix with integer entries $A=[a_{ij}]\in \Z^{m\times n}$, and $\x\in\R^n$.

 For a nonzero vector $\x=(x_1,\ldots,x_n)\trans\in \R^n$ we define the \emph{relative magnitude} of $\x$ as
 \begin{equation}\label{defcondnum}
 \omega(\x)=\max\left\{\frac{|x_j|}{|x_i|}, \x\in \R^n\setminus\{\0\}, x_i\ne 0\right\}.
 \end{equation}

 It is easy to check that if $\x$ has no zero coordinates then the relative magnitude of a vector $\x$ coincides with the classic condition number of the diagonal matrix that has the entries of $\x$ in its diagonal.

 Denote by $N(A)$ the nullspace of $A$ and let $N'(A):=N(A)\setminus \{0\}$. 
 
 The \emph{solution relative magnitude} of $A$ is given by
 \begin{equation}\label{solconnum}
 \omega(A)=\inf\{\omega(\x),\;\x\in N'(A) \}.
 \end{equation}
 We agree that if $A\x=\0$ has the unique solution $\x=\0$ then $\omega(A)=0$.

 The aim of this paper is to establish a sharp upper bound on $\omega(A)$ in terms of $\|A\|_{\infty}:=\max_i \sum_{j=1}^n |a_{ij}|$. Namely, we have the following result.
 \begin{theo} Consider a nonzero matrix $A\in \Z^{m \times n}$. If $\|A\|_\infty=1,2$, then $\omega(A)=0$ or $\omega(A)=1$. For $\|A\|_\infty \geq 3$, we have the following sharp upper bound for the solution relative magnitude of $A$:
 \label{firsttheo}
 \begin{equation}\label{mainineq}
 \omega(A)\le (\|A\|_{\infty}-1)^{\rank A}.
 \end{equation}
\end{theo}

 It is quite easy to see the sharpness of \eqref{mainineq}. Let $A\x=\0$ be the system of $n-1$ homogeneous equations  $k x_i -x_{i+1}=0$ for $i=1,\ldots,n-1$ for a given $k\in\N$.   Then $\|A\|_{\infty}=k+1, \rank A=n-1$ and $\omega(A)=\omega(\x)=\frac{x_n}{x_1}=k^{n-1}$ for
 any $\x\ne\0$ in the null space of $A$.

 The cases $\|A\|_{\infty}=1,2$ are simple (cf.\ Proposition \ref{casesnormA12}).  The case $\|A\|_{\infty}\ge 3$ is deduced from the following result.
 \begin{theo}\label{maintheo}  Fix an integer $k \ge 2$.
 Consider a linear system in which all equations are of one of the following two types:
 \begin{eqnarray}
 && x_i=\pm 1 \label{first_type}\\
 && \pm x_{i_1} \pm x_{i_2} \pm \ldots \pm x_{i_{l+1}} =0,\label{genformeq}
 \end{eqnarray}
 where $l$ is a non-negative integer satisfying $l\le k$.  (The integer $l$ may depend on the equation.)
 Assume that this system is solvable. Then there exists a rational solution such that $|x_j|\leq k^{n-1}$ for each variable $x_j$,  with $n$ being the rank of the system. The above bound is sharp.
 \end{theo}
 For the system $x_1=1$ and $k x_i -x_{i+1}=0$ for $i=1,\ldots,n-1$ our theorem is sharp.
 Our main tool is the Hadamard-Fischer determinant inequalities combined with graph theoretical arguments.

 The case $k=2$ of Theorem \ref{maintheo} proves one of the conjectures of A.\ Tyszka presented in \cite{Tys09,Tys10}.
 See \S2.\medskip

 We now survey briefly the contents of our paper. In \S2 we discuss some properties of $\omega(A)$ and show that Theorem \ref{maintheo} implies Theorem \ref{firsttheo} and one the above mentioned
 Tyska's conjectures.  In \S3 we lay the ground for the proof of Theorem \ref{maintheo}.   We reduce the system
 \eqref{first_type}--\eqref{genformeq}
 into a system of the same type of equations with smaller number of variables satisfying the following conditions.
 First, the system has a unique solution; second, the system \eqref{first_type} is $x_1=1$; third, no variable is equal to zero, fourth, $x_i\ne \pm x_j$  for $i\ne j$.  These conditions allow us to split the new system of equations \eqref{genformeq} to a finite number of maximal chain equations of the form $kx_{i_{j+1}}=\pm x_{i_j}$ for $j=1,\ldots,t$. We show that no two maximal chains have a common variable.
 In \S4 we estimate the determinants of certain tridiagonal matrices related to a maximal chain, and the Euclidean norm of the coefficients of
 each equation in \eqref{genformeq} which does not appear in any maximal chain.  In \S5 we conclude the proof of Theorem \ref{maintheo}.

 \section{On Relative Magnitudes}
 For $m\in\N$ denote $[m]:=\{1,\ldots,m\}$.
 Let $\x=(x_1,\ldots,x_n)\trans\in\R^n\setminus\{\0\}$.  Define the {\em relative magnitude}\/ of $\x$ by \eqref{defcondnum}.  Let
 \[\phi(\x)=\min\{|x_i|, x_i\ne 0\}, \quad \Phi(\x)=\max\{|x_i|, i\in [n]\}.\]
 So $\omega(\x)=\frac{\Phi(\x)}{\phi(\x)}$.
 Let $A\in\R^{m\times n}$.  If $\rank A\le n-1$
 we define the solution relative magnitude of $A$ by \eqref{solconnum}, otherwise we let $\omega(A)=0$.
 The following result shows that the infimum in \eqref{solconnum} is attained.
 \begin{prop}\label{achsolcondn}
 Let $A\in \R^{m\times n}$ and assume that $\rank A\le n-1$.
 Then
 \begin{equation}
 \label{achsolcondn1}
 \omega(A)=\min \{ \omega(\x),\x\in N'(A)\} \geq 1.
 \end{equation}
  So, we can write
 $$\omega(A) = \min_{\x\in N'(A)} \max_{x_i\neq 0} \frac{|x_j|}{|x_i|}.$$
 \end{prop}
 \proof Under these conditions, we must have a nonzero solution $\z$ of the system $A\x=\0$. If $z_i$ is a nonzero coordinate of $\z$, then $\omega(\z)\geq |z_i|/|z_i| =1$. Since this holds for any nonzero solution, $\omega(A)\geq 1$.

 Now let
 $$\Sigma(A):=\left\{\x=(x_1,\ldots,x_n)\trans\in\R^n, \x\in N(A), \sum_{i=1}^n |x_i|=1\right\}.$$
 We clearly have $\omega(A)=\inf\{\omega(\x), \x\in \Sigma(A)\}$, and $\Sigma(A)$ is a compact set, but the function $\omega$ is in general not continuous on $\Sigma(A)$. For instance, if $A=[1\ 0\ 0]$ then $(0,1,0)\trans,(0,1-t,t)\trans\in \Sigma(A)$, $(0,t,1-t)\trans\in \Sigma(A)$ if and only if $t\in [0,1]$, as $|t|\le 1$ and $|1-t|\le 1$.
 
For $t\in (0,\frac{1}{2})$, we have $\omega((0,1,0)\trans)=1$, $\omega((0,1-t,t)\trans)=(1-t)/t$, $\omega((0,t,1-t)\trans))=\frac{t}{1-t}$ 
 
To circumvent this problem, we start by noticing that $\omega(A)=\inf\{\omega(\x), \x\in \Sigma(A), \omega(\x)\le 2\omega(A)\}$. Assume that $\x\in \Sigma(A), \omega(\x)\le 2\omega(A)$.
 As $1\le (n-1)\Phi(\x)+\phi(\x)$ we deduce that
 \[\frac{1}{\phi(\x)}\le (n-1)\omega(\x)+1\le 2(n-1)\omega(A)+1,\]
and therefore
\[\phi(\x)\ge \alpha(A):=\frac{1}{2(n-1)\omega(A)+1}.\]
 Let $\Sigma_1(A):=\{\x\in\Sigma(A), \phi(\x)\ge \alpha(A)\}$.  Clearly, $\Sigma_1(A)$ is a compact set, and both $\Phi(\x)$ and $\phi(\x)$ are continuous on $\Sigma_1(A)$. Hence $\omega(\x)$ is a continuous function on $\Sigma_1(A)$ and
 $$\omega(A)=\inf\{\x\in\Sigma_1(A)\} = \min\{\x\in\Sigma_1(A)\}.$$
 Hence, condition \eqref{achsolcondn1} holds.\qed\medskip

Define the {\it support}\,\,Êof a vector $\x$ as the set $I\subseteq [n]$ such that $x_i\neq 0 \eq i\in I$. Following \cite{BFP}, given a subspace $W\leq \R^n$, we say that a nonzero vector is {\em elementary}\ (in $W$) if its support is minimal among all supports of nonzero vectors in $W$. In other words $\x\in W$ is elementary if for $\y\in W$, $\y\neq \0$ and $\supp \y \subseteq \supp \x$, we have $\supp \y = \supp \x$. 
It is known that any subspace has a basis formed by elementary vectors (see \cite{BFP} and \cite[p.\ 528]{CP} for an algorithm to find such a basis). 

 \begin{prop}\label{casesnormA12}  Let $A\in \Z^{m\times n}$ and assume that $\|A\|_{\infty}=1,2$. Then, for every elementary vector $\x\in N(A)$, we have $\omega(\x)=1$. Therefore $\omega(A)=0$ or 1, as stated in Theorem \ref{firsttheo}.
 \end{prop}
 \proof  Suppose first that $\rank A =n$.  Then $\omega(A)=0$.
 Assume now that $\rank A\le n-1$.  Suppose first that $\|A\|_{\infty}=1$.  Then each nontrivial equation of $A\z=\0$ is given by $x_i=0$ for some
 $i\in [n]$.  Hence the general solution of $A\z=\0$ is of the following form.  The set of free variables is a nonempty strict subset $S$
 of $[n]$ and all other dependent variables equal to zero.
 In particular, all elementary vectors have only one nonzero coordinate: $A\e_j=\0$ for $j\in S$, hence $\omega(\x)=\omega(A)=1$.

 Assume now that $\|A\|_{\infty}=2$.  Then the nontrivial equations of $A\z=0$ are either of the form $x_i=0$ or $x_i=\pm x_j$.
 If $A\e_j=\0$ for some $j\in [n]$ then $\omega(\e_j)=\omega(A)=1$.
 Assume that $A\e_j\ne\0$ for each $j\in[n]$.
 It is easy to see that the general solution of this system has the following general form. One can partition the set $[n]$ into
 $$[n]= S\, \dot\cup\, T_1\, \dot\cup\, \ldots\, \dot\cup\, T_l,$$
 so that $x_i=0$ if $i\in S$, each $T_p$ contains at least two elements, for $p\in [l]$ and, for each pair of distinct indices $i,j\in T_p$, we have $|x_i|=|x_j|$.
 The value of all $|x_i|, i\in T_p$, can be prescribed arbitrarily for each $p\in[l]$; in particular, elementary vectors $\x$ satisfy $x_i=0$ for $i\not\in T_p$ and $|x_j|\neq 0$ for $j\in T_p$, for some $p\in [l]$.  Hence $\omega(\x)=\omega(A)=1$.  \qed\medskip

The following two results prove that, once Theorem \ref{maintheo} is proved, then Theorem \ref{firsttheo} holds also for $\|A\|_\infty \geq 3$.

\begin{prop}
Let $A\in \Z^{m\times n}$ with $\rank A\leq n-1$ and $\|A\|_\infty \geq 3$, and let $\x\neq \0$ be an elementary vector of $N(A)$, with support $I$. Then $|I|-1\leq \rank A$ and 
$$\omega(\x)\leq (\|A\|_\infty-1)^{|I|-1}.$$

In particular, $\omega(A)\leq (\|A\|_\infty-1)^{\rank A}$, as stated in Theorem \ref{firsttheo}.
\end{prop}
\proof If $|I|=1$, then the result holds. Now suppose $|I|\geq 1$ and let $\y$ be the vector of $\R^{|I|}$ formed by the nonzero coordinates of $\x$, $\omega(\x)=\omega(\y)$. Let $B$ be the matrix obtained from $A$ by selecting the columns with indices in $I$. Then $B\y=\0$ and we must have $\rank B=|I|-1$, otherwise we would be able to get a solution of $A\z=\0$ with more zero coordinates (by setting one of the free variables to zero), contradicting the fact that $\x$ is elementary. Clearly $|I|-1\leq \rank A$. Set $k+1:=\|A\|_\infty\geq \|B\|_\infty$, $k\geq 2$.

Since the null space of $B$ is spanned by one vector $\w\ne \0$ we deduce that $\omega(\w)=\omega(\y)$. Without loss of generality we may assume that $\phi(\w)=1$, with $w_j=1$. If we now consider the system $B\z=\0$ along with the equation $z_j=1$, its only solution will be $\y$ and the system is of the type described in Theorem \ref{maintheo}, with rank $|I|$. Once this result is proved, we get that $\omega(\y)=\Phi(\y) \leq k^{|I|-1}$, and hence
$$\omega(\x)=\omega(\y)\leq (\|A\|_\infty-1)^{|I|-1}.$$
This proves the result.\qed\medskip

This previous result allows for an improvement of Theorem \ref{maintheo}.

\begin{corol}
Consider a nonzero matrix $A\in \Z^{m \times n}$ with $\|A\|_\infty \geq 2$. Let $t>0$ be the least number of nonzero coordinates in a nonzero vector of $N(A)$. Then $t-1\leq \rank A$ and
$$\omega(A)\leq (\|A\|_\infty-1)^{t-1}\leq (\|A\|_\infty-1)^{\rank A}.$$
\end{corol}

 One of the conjectures of A.\ Tyszka, presented in \cite{Tys09,Tys10}, is:
\begin{con}\label{mthm}
Let $n\ge 2$. Assume that we have a solvable linear system of $n$ equations of the following two types:
\begin{equation}\label{ijkcon}
x_i+x_j=x_k \textrm{\quad and\quad } x_l=1.
\end{equation}
Then the above system has a rational solution with $|x_i|\leq2^{n-1}$ for all $1\leq i\leq n$.  This result is sharp.
\end{con}

Some partial results about this conjecture were obtained by Tyszka \cite{Tys09}, with upper bound $\sqrt{5}^{\, n-1}$, and by Cipu \cite{Ci}, with upper bound $2^n$.  Clearly, Theorem \ref{maintheo} for $k=2$ yields Conjecture \ref{mthm}.
The rest of this paper is devoted to proving Theorem \ref{maintheo}.

 \section{Simplification}
Assume that our given system of equations has variables $x_i, i\in[m]$. We'll show that it is enough to prove Theorem \ref{maintheo} after the following simplifications.\medskip

 1. Observe first that if all equations are of type \eqref{genformeq}, we have a homogeneous system, for which there is a zero solution. Now we assume we have a nonempty set of equations $x_i=\pm 1$ for $i\in S$, and $|S|\ge 2$. We can replace this set by
 one equation $x_i=\pm 1$ and equations $x_j \pm x_i=0$ for $j\in S\setminus\{i\}$.\smallskip

 2. If the system is indeterminate, let $S\subset \N$ be a set of free variables, and for every $x_j,j\in S$, take $x_j=0$. We get new equations of the type \eqref{genformeq} with variables indexed in $[m]\setminus S$. Out of those, we select any $m-|S|-1$ linearly independent equations in $m-|S|$ variables.
 Rename these variables so that the first equation is $x_1=\pm 1$ and the other variables are $x_j, j=2,\ldots,m'=m-|S|$.
 Hence the new system in $m'$ variables has a unique solution $\x=(x_1,\ldots,x_{m'})\trans$.\smallskip

 3. Let $T\subset [m']$ be the set of all indices $j$ for which $x_j=0$.
 As above we replace the above system by a system with variables indexed in $[m']\setminus T$, where the system has a unique solution and each 
 $x_j\ne 0$ for $j\in [m']\setminus T'$.  Again we rename our variables so now we have $n$ variables $x_1,\ldots,x_n$, where the first equation
 is $x_1=\pm 1$ and all other $n-1$ equations are of the form \eqref{genformeq}.  This system has a unique solution and each $x_j\ne 0$.
 Note that the rank of this system will be less than or equal to the rank of the original system.\smallskip

 4.  Suppose that for this unique solution $\x$ we have $x_i=\pm x_j$ for $i< j$.  Then we eliminate $x_j$ in all equations \eqref{genformeq} by substitution $x_j=\pm x_i$.  Note that we still have a unique solution of the system with one variable less.

 We repeat this process until we have a system $x_1=\pm 1$, all other equations are of the form \eqref{genformeq}, such that this system
 has a unique solution and $x_i\ne \pm x_j$ for $i\ne j$.  So again we can assume that we have a system in which all equations are of the form \eqref{first_type} or \eqref{genformeq}.\smallskip

 5. Next we cancel terms in the left-hand side of \eqref{genformeq}.  I.e. if we have a term $x_j$ and a term $-x_j$ in the left-hand side of \eqref{genformeq}, then we replace the two terms by $0$.  Thus we can assume that the left-hand side of \eqref{genformeq} does not contain the
 same variable $x_j$ with opposite signs.  Since each solution $x_j\ne 0$, and the system has a unique solution, it means that after the cancelations each of $n-1$ equations of the form \eqref{genformeq}  contains at least two variables.

 To summarize, we reduced our problem to the following system of $n$ equations with $n$ variables $x_1,\ldots,x_n$ one of them of the form $x_1=\pm 1$
 and the other $n-1$ equations of the form
 \begin{eqnarray}\label{genformeq1}
 &&\sum_{j=1}^{l+1} \pm c_{i_j}x_{i_j}=0, \quad 1\le i_1<\ldots<i_{l+1}\le n,\\
 && c_1,\ldots,c_{l+1} \in \N, 1\le l\le k, \textrm{ and }c_1+\ldots+c_{l+1}\le k+1.\label{genformeq2}
 \end{eqnarray}
 Furthermore the above system has a unique solution $\x=(x_1,\ldots,x_n)\trans$ with the following properties:
 $x_i\ne 0$ for each $i$ and $x_i\ne \pm x_j$ for each $i\ne j$.
 We note that not all such systems are being considered. For instance, we
 are not allowing the case $l=1$ and $c_1=c_2$, which would mean that $x_{i_1}=\pm x_{i_2}$.\smallskip

We note that these reductions do not change the maximum value of $|x_p|$, but may reduce the rank of the original system.  Hence it is enough to prove Theorem
 \ref{maintheo} for the system we have obtained. \smallskip

 6. Among the homogeneous equations, we now consider those of type $kx_i=\pm x_j$, with $i\neq j$. Call these {\em equations of type 2}.  The equation $x_1=\pm 1$ is called of {\em type 1} and all others will be of {\em type 3}. By replacing $x_1$ by $\pm x_1$ we can assume that the equation of type 1 is $x_1=1$.
 We claim that we can organize equations of type 2 in maximal chains of the form
 $$kx_{i_2}=\pm x_{i_1},\quad kx_{i_3}=\pm x_{i_2},\quad \ldots \quad kx_{i_t}=\pm x_{i_{t-1}},$$
 which we denote as the chain $i_1\to i_2\to\ldots\to i_t$.   So $i_1$ and $i_t$ are the \emph{head} and the \emph{tail} of the chain, respectively.
 We now claim that no variable appears in more than one maximal chain. \medskip

 Assume to the contrary that two maximal chains intersect at $i$.  Since these two chains are maximal $i$ is not the tail of one chain
 and the head of another chain.
 Thus we only have the following two cases which are impossible in view of our assumptions.

 \begin{itemize}
 \item $kx_i=\pm x_j$ and $kx_i= \pm x_p$. Hence $x_j=\pm x_p$.  This violates the condition that $x_j\ne \pm x_p$
 \item $kx_j=\pm x_i$ and $kx_p= \pm x_i$. Hence $x_j=\pm x_p$.  This violates the condition that $x_j\ne \pm x_p$
 \end{itemize}

 We are now left with a system of equations with a unique solution, in which no chains intersect.\medskip

 Reorder the variables so that in each chain we have
 $$kx_{i_1}=\pm x_{i_1+1},\quad kx_{i_1+1}=\pm x_{i_1+2},\quad \ldots \quad kx_{i_1+t-1}=\pm x_{i_1+t}.$$

 After this, when we group all equations belonging to the same chain, the system matrix will have a block of rows that looks like this:

 $$\begin{bmatrix}
 k & \pm 1 \\
 & k & \pm 1\\
 && k & \pm 1\\
 &&& \ddots & \ddots \\
 &&& & k & \pm 1
 \end{bmatrix}$$

 This block is of type $t\times (t+1)$, that is to say, the number of columns is the number of rows plus one.\medskip

 Before these blocks of equations, we put first the equation of type 1. Since we have renumbered the variables, it is no longer necessary that the variable in the equation of type 1 has subindex 1, the equation now is $x_s=1$ for some index $s$. Moreover, it is possible that the variable $x_s$ belongs to some maximal chain. Then we list all equations of type 3. We note that if we have $r$ chains, we must have at least $r-1$ equations of type 3, so that we get a square matrix. Call the resulting matrix $A$.\medskip

 We now estimate the magnitude of the (only) solution of our system using Cramer's rule. Let $A_i$ be the matrix obtained by replacing column $i$ of $A$ by $e_1$, the column of independent terms. Since $|\det A|\geq 1$, we need to establish that $|\det A_i|\leq k^{n-1}$, where $n$ is now the number of variables and equations of the system, and the size of the matrix $A$.\medskip

 Let $U_i$ be the $(n-1)\times (n-1)$ submatrix obtained from $A_i$ by deleting the first row and the $i$-th column.   We claim that $|\det A_i|=|\det U_i|$
 for $i=1,\ldots,n$.  Indeed, since the first equation in our system is $x_s=1$ then expanding $\det A_s$ by the first row we obtain that $\det A_s=\det U_s$. Now take $i\neq s$.  Subtract the column $i$ in $A_i$ from the column $s$ to obtain the matrix $V_i$.  Expand the determinant of $V_i$ by the first row
 to deduce $|\det A_i|=|\det V_i|=|\det U_i|$.

 \section{Principal minors of $U_i U_i\trans$}

 To obtain a bound for $|\det U_i|$, we will consider $\det W_i$, where $W_i=\emph{}U_i U_i\trans$, which is majorized by a product of its principal minors according to the Hadamard-Fischer inequality \cite[Th.\ 7.8.3, p.\ 478]{HJ}. In $W_i$, we consider the principal minors determined by the row blocks we have defined.\medskip

 Assume that the block of rows corresponding to a chain with $t$ rows is uncut by the column $i$; i.e., the variable $x_i$ is not participating in this
 chain.  Then for $t\ge 2$ the corresponding principal submatrix of $W_i$ is
 a $t\times t$ tridiagonal symmetric matrix:
 $$B_t:=\begin{bmatrix}
 k^2+1 & \pm k \\
 \pm k & k^2+1 & \pm k\\
 & \pm k & k^2+1 \\
 &&  & \ddots & \pm k \\
 &&& \pm k & k^2+1
 \end{bmatrix} \in \Z^{t\times t}.$$
 For $t=1$ we have $B_1=[k^2+1]$.

 We now consider the case where the deleted column $i$ in A
 does cut through a block of rows in a chain.   Then the corresponding principal submatrix of $W_i$ is the direct sum of two blocks:
 \begin{eqnarray*}
 & C_p:=\begin{bmatrix}
 k^2+1 & \pm k \\
 \pm k & \ddots  \\
 &  & k^2+1 & \pm k\\
 && \pm k & k^2+1 & \pm k\\
 &&& \pm k & k^2
 \end{bmatrix}\in \Z^{p\times p},\\
 & D_q:=\begin{bmatrix}
 1 & \pm k \\
 \pm k & k^2+1 & \pm k\\
 & \pm k & k^2+1 \\
 &&  & \ddots & \pm k \\
 &&& \pm k & k^2+1
 \end{bmatrix} \in \Z^{q\times q},
 \end{eqnarray*}
 where $p,q\ge 2$ and $p+q=t$.  Note that $C_1=k^2, D_1=1$.
 It is possible that $p$ or $q$ is zero, i.e., we have only one block instead of a direct sum of two blocks.
 We define $\det C_0=\det D_0=1$.

 \begin{lemma}\label{formdetBCDt}
 For $t\in\N$ and a real number $k$ the following equalities hold.
 \begin{eqnarray}\label{formdetBt}
 && \det B_t=\frac{(k^2)^{t+1}-1}{k^2-1},\\
 && \det C_t=k^{2t}, \label{formdetCt}\\
 && \det D_t=1.\label{formdetDt}
 \end{eqnarray}
 ($\det B_t=(t+1)$ for $k=\pm 1$.)
 \end{lemma}
 \proof  It is enough to consider the case $k\ne \pm 1$.  For $t=1,2$ \eqref{formdetBt}, \eqref{formdetCt} and  \eqref{formdetDt} clearly hold.
 Assume that $t\ge 3$.
 Using the Laplace expansion by the first row of $B_t,D_t$ and by the last row of $C_t$ we obtain
 \begin{eqnarray*}
 \det B_{t}=(k^2+1)\det B_{t-1}-k^2 \det B_{t-2},\\
 \det C_t=k^2\det B_{t-1}-k^2\det B_{t-2},\\
 \det D_t=\det B_{t-1}-k^2\det B_{t-2}.
 \end{eqnarray*}
 Consider first the recurrence equations for $\det B_t$, $t\ge 3$.
 The roots of the characteristic polynomial of this recurrent system are 1 and $k^2$.
 So $\det B_t=a_1(k^2)^t +a_0$.  Since  \eqref{formdetBt} is of this form and it holds for $t=1,2$,  \eqref{formdetBt}
 holds for all $t\in\N$.

 Substitute the expression \eqref{formdetBt} in the expressions for $\det C_t$ and $\det D_t$ to deduce \eqref{formdetCt} and \eqref{formdetDt}
 respectively.  \qed\medskip

 We now discuss the equations of type 3.
 \begin{lemma}\label{upboundtype3}  Let $2\le k\in\N$.
 Consider the equation \eqref{genformeq1} of type 3, i.e. $l\in [2,k]$ and which is not of the form $kx_{i}\pm x_j=0, i\ne j$.
 Then the $\ell_2$ norm of the coefficient vector $\|(\pm c_{i_1},\ldots,\pm c_{i_{l+1}})\|_2=\sqrt{\sum_{j=1}^{l+1} c_{i_j}^2}$  is at most $\min(\sqrt{k^2-1},\sqrt{(k-1)^2+4})\leq \sqrt{k^2-1}$.

 Furthermore, if we delete a variable $x_{i_j}$ from \eqref{genformeq1} for some integer $j\in [1,l+1]$ then
 the $\ell_2$ norm of the coefficient vector of this new equation is at most $\sqrt{(k-1)^2+1}$.

 In particular, consider the diagonal element $w_{tt}$ of $U_iU_i\trans$ corresponding to the equation \eqref{genformeq1} of type 3.
 If the variable $x_i$ does not appear in this equation, $w_{tt}\le k^2-1$.  If $x_i$ appears in this equation then $w_{tt}\le (k-1)^2 +1$.
 \end{lemma}
 \proof  Let $x,y\ge 0$ and assume that $z=x+y$.  Suppose that $a\in [0, \frac{z}{2}]$ and $x,y\ge a$.  It is straightforward to show that $x^2+y^2\le (z-a)^2 +a^2$, with $x^2+y^2=(z-a)^2+a^2$ if and only if either $x=z-a,y=a$ or $x=a,y=z-a$.

 Suppose first that $k=2$.  Then the equation \eqref{genformeq1} of type 3 is either $\pm x_{i_1}\pm x_{i_2}=0$ or $\pm x_{i_1}\pm x_{i_2}\emph{}\pm x_{i_3}=0$.
 Then the $\ell_2$ norm of the coefficient vector of this system is at most $\sqrt{3}=\sqrt{k^2-1}< \sqrt{(k-1)^2+4}$.

 If we delete one of the variables in this equation, then the $\ell_2$ norm of the coefficient vector of this new equation is at most $\sqrt{2}=\sqrt{(k-1)^2+1}$.

 Suppose now that $k\ge 3$.  Assume first that $l\ge 2$.  Rename the indices so that $c_{i_l},c_{i_{l+1}}\in [1,k-2]\cap \N$.  Then $c_{i_l}^2 + c_{i_{l+1}}^2 < (c_{i_{l}}+c_{i_{l+1}})^2 +0^2$.  Hence
 $$c:=\sum_{j=1}^{l+1} c_{i_j}^2< (c_{i_{l}}+c_{i_{l+1}})^2 +\sum_{j=1}^{l-1} c_{i_j}^2.$$

 It now follows that the maximal value of $c$ is achieved when the equation has only two variables, which correspond
 to the value $l=1$.  So the maximum value of $c$ is achieved when the coefficient of one variable is $\pm (k-1)$ and the coefficient of the other variable is $\pm 2$.

 Assume now that in \eqref{genformeq1} $l=1$, i.e. $\pm c_{i_1}x_{i_1} \pm c_{i_2} x_{i_2}=0$, where $c_{i_1},c_{i_2}\in
 [1,k-1]\cap \N$ and $z=c_{i_1}+c_{i_2}\le k+1$.  Using our observation in the beginning of the proof of this lemma
 it is straightforward to show that $c_{i_1}^2 +c_{i_2}^2\le (k-1)^2 +2^2$ (equality holds if and only if $c_{i_1}=k-1, c_{i_2}=2$ or $c_{i_1}=2, c_{i_2}=k-1$). For $k>2$, $(k-1)^2+4\leq k^2-1$.

 Assume now that we remove one variable from \eqref{genformeq1}.  By renaming the variables we may assume that it is the variable $x_{i_{l+1}}$.
 Thus we need to find an upper bound for $\sum_{j=1}^l c_{i_j}^2$, where $c_{i_1},\ldots,c_{i_{l}}\in [1,k-1]\cap \N$ and $\sum_{j=1}^l c_{i_j}\le k$.
 For $l=1$ the upper bound is $(k-1)^2$.  For $l\ge 2$, from the above arguments, we deduce that this upper bound is $(k-1)^2 +1$.
 Equality holds if and only if $l=2$, $\{c_{i_1},c_{i_2},c_{i_3}\}=\{k-1,1,1\}$ and the coefficient of the deleted variable $x_{i_j}$ is $1$.  \qed

 \section{Proof of Theorem \ref{maintheo}}
 Assume that after the reductions described in \S2 we have a system of $n$ linearly independent equations
 of the following forms.  The first equation is
 of the form $x_s=1$. The other $n-1$ equations are of the form \eqref{genformeq1}, which are of types 2 and 3.

 The equations of type 2 are organized in groups.  Each group of type 2 equations is a system of $t$ equations in $t+1$ variables of the form $kx_{i_j}=\pm x_{i_{j+1}}, j=1,\ldots,t$, which is called a chain of length $t+1$.   No two distinct chains have a common variable.

 The number of chains is $r\ge 0$.  The equations of the type 3 are of the form \eqref{genformeq1}, where $c_{i_j}\in \N, c_{i_j}\le k-1$ for $j=1,\ldots,l+1\le k+1$.  Since we have a unique solution to our system, we must have at least $r-1$ equations of type $3$.

 Hence our system of equations is of the form $A\x=\e_1, A\in \Z^{n\times n}, |\det A|\ge 1$.  Recall that we denoted by $U_i$ the submatrix of $A$ obtained by deleting
 the first row and the column $i$ for $i=1,\ldots,n$.  Then Cramer's rule yields that $|x_i|\le |\det U_i|$.\medskip

 Also, recall that in one of our reductions we reduced the number of variables whenever we had the equality $x_i=\pm x_j$ for $i\ne j$ or $x_i=0$. These reductions do not change the maximum value of $|x_p|$, but may reduce the rank of the original system.  Hence it is enough to prove Theorem
 \ref{maintheo} in the above form.\medskip

 We prove Theorem \ref{maintheo} by showing that
 \begin{equation}\label{fundin}
 x_i^2\le \det W_i=\det U_iU_i\trans\leq k^{2(n-1)}, i=1,\ldots,n,
 \end{equation}
  using the Hadamard-Fischer inequality and Lemmas \ref{formdetBCDt} and \ref{upboundtype3}.
  This is done by considering two cases.
 \medskip

 {\bf Case 1}. Column $i$ cuts through a chain block of length $t+1$.
 This block yields a principal minor in $W_i$, which will be less than or equal to $k^{2t}=\det C_t\det D_0\ge \det C_p\det D_q$, $p+q=t$, by Lemma \ref{formdetBCDt}.

 Assume first that we have only one chain.  If
 this chain is of length $n$ we just showed that \eqref{fundin} holds.  If $t<n-1$ then Lemma \ref{upboundtype3} and the Hadamard-Fischer inequality yield
 $\det W_i\le k^{2t} (k^2-1)^{n-1-t}<k^{2(n-1)}$.

 Assume now that we have $r\ge 2$ chains.
 The lengths of the chains are $t_1+1,\ldots, t_{r-1}+1, t_r+1=t+1$.
 Lemma \ref{formdetBCDt} yields that the principal minor of $W_i$ corresponding to the chain of length $t_j+1$ is bounded above by
 $\frac{k^{2(t_j+1)}-1}{k^2-1}$.  Therefore the product of the $r-1$ principal minors in $W_i$ corresponding to chains $1,\ldots,r-1$ is bounded above by
 \[\prod_{j=1}^{r-1} \frac{k^{2(t_j+1)}-1}{k^2-1} < \frac{k^{2\sum_{j=1}^{r-1} (t_j+1)}}{(k^2-1)^{r-1}}.\]
 The number of the equations of type $3$ is $n-1-\sum_{j=1}^{r} t_j\ge r-1$.  In view of Lemma \ref{formdetBCDt} the product of the diagonal
 entries in $W_i$ corresponding to the equations of type $3$ is bounded above by $(k^2-1)^{n-1-\sum_{j=1}^{r} t_j}$.
 Combine all the above inequalities to deduce that
 \begin{eqnarray*}
 \det W_i & < & k^{2t_r}\frac{k^{2\sum_{j=1}^{r-1} (t_j+1)}}{(k^2-1)^{r-1}} (k^2-1)^{n-1-\sum_{j=1}^{r} t_j} \\
 & \leq & k^{2(t_r+\sum_{j=1}^{r-1}(t_j+1))}(k^2-1)^{n-1-(r-1)-\sum_{j=1}^{r} t_j}\\
 & \le & k^{2(n-1)}.
 \end{eqnarray*}

 {\bf Case 2}. Column $i$ cuts through a row of type 3 but not through any chain.

 If we have only one chain of length $t+1$, then
 $$\det W_i \leq \frac{k^{2(t+1)}-1}{k^2-1} (k^2-1)^{n-t-1} < k^{2(t+1)} k^{2(n-t-2)} = k^{2(n-1)}.$$
 If there is no chain, the analysis is similar.

 Assume now that we have $r\ge 2$ chains of lengths $t_1+1,\ldots,t_r+1$.  Hence we have $n-1-\sum_{j=1}^r t_j\ge r-1$  equations of type 3.

 One equation of type 3 discussed above must contain the variable $x_i$.
 Since $(k-1)^2+1\le k^2 -2$ for $k\ge 2$, we can conclude, using Lemma \ref{upboundtype3}, that the diagonal entry of $W_i$ corresponding to this equation does not exceed $k^2-2$.
% Lemma \ref{upboundtype3} also states that the diagonal entries of $W_i$ corresponding to all other equations of type 3 do not exceed $k^2-1$.
 Hence $\det W_i$ is less that or equal to
 \[
 \frac{k^{2(t_1+1)}-1}{k^2-1} \frac{k^{2(t_2+1)}-1}{k^2-1}(k^2-2)(\prod_{j=3}^r  \frac{k^{2(t_j+1)}-1}{k^2-1})(k^2-1)^{n-2-\sum_{j=1}^r t_j}.\]
 Recall that  $n-2-\sum_{j=1}^r t_j\ge r-2$.  Hence as in the Case 1 we get
 \begin{eqnarray*}
 &&(\prod_{j=3}^r  \frac{k^{2(t_j+1)}-1}{k^2-1})(k^2-1)^{n-2-\sum_{j=1}^r t_j}\\
 &&\leq k^{2\sum_{j=3}^r (t_j+1)}(k^2-1)^{n-2-\sum_{j=1}^r t_j-(r-2)}\le k^{2(n-2-(t_1+t_2))}.
 \end{eqnarray*}

 As $(k^2-1)^2>k^2(k^2-2)$ it follows that
 \begin{eqnarray*}
 \frac{k^{2(t_1+1)}-1}{k^2-1} \frac{k^{2(t_2+1)}-1}{k^2-1}(k^2-2) & < & \frac{(k^{2(t_1+1)}-1)(k^{2(t_2+1)}-1)}{k^2(k^2-2)} (k^2-2)\\
 & < & \frac{k^{2(t_1+1)}k^{2(t_2+1)}}{k^2}=k^{2(t_1+t_2)+2}.
 \end{eqnarray*}

 Combine the above inequalities to deduce $\det W_i\le k^{2(n-1)}$.
 The proof of Theorem \ref{maintheo} is concluded.  \qed\medskip
 
We thank the referees for their useful remarks.

 \bibliographystyle{plain}

\end{document}